\documentclass[10pt]{article}
\usepackage{latexsym}
\usepackage{amsfonts}

\input amssym.def
\input amssym.tex

\newcommand{\ZZ}{{\mathbb Z }}
\newcommand{\QQ}{{\mathbb Q }}
\newcommand{\RR}{{\mathbb R }}
\newcommand{\CC}{{\mathbb C }}
\newcommand{\PP}{{\mathbb P }}

\newcommand{\renorm}{{ \setcounter{equation}{0} }}

\title{Cohomology of Moduli Spaces}
\author{Frances Kirwan\\Oxford University, UK}
\date{}
\begin{document}
\maketitle

Moduli spaces arise in classification problems in algebraic
geometry (and other areas of geometry) when, as is typically
the case, there are not enough discrete invariants to classify
objects up to isomorphism. In the case of nonsingular complex
projective curves (or compact Riemann surfaces) the genus $g$ is
a discrete invariant which classifies the curve regarded as
a topological surface,
but does not determine its complex structure when $g>0$. For each
$g\geq 0$ there is a moduli space ${\mathcal M}_g$ 
whose points correspond bijectively to isomorphism classes of
nonsingular complex projective curves of genus $g$, and whose
geometric structure reflects the way such curves can vary in
families depending on parameters.
The topology of these moduli
spaces ${\mathcal M}_g$ and their compactifications has been 
studied for several decades, and 
 important progress has been made recently on some long-standing questions concerning their cohomology. 

In his fundamental paper \cite{Mumford} Mumford considered some
tautological cohomological classes $\kappa_j \in H^{2j}({\mathcal M}_g)$
for $j=1,2,\ldots$ which extend naturally to the
Deligne-Mumford compactification $\overline{{\mathcal M}}_g$. 
Much work on the cohomology of ${\mathcal M}_g$ 
 has concentrated on its tautological
ring, which is the subalgebra of its rational cohomology ring
 (or of its Chow ring)
generated by these tautological classes. 

One reason for the importance
of the tautological ring of ${\mathcal M}_g$ is its relationship
with the stable cohomology ring $H^*({\mathcal M}_\infty;\QQ)$. It was proved by Harer \cite{Harer}
that the cohomology $H^k({\mathcal M}_g;\QQ)$ of ${\mathcal M}_g$ in
degree $k$ is independent of the genus $g$ when $g \gg k$, making it possible to define
$H^k({\mathcal M}_\infty;\QQ)$ as $H^k({\mathcal M}_g;\QQ)$ for suitably
large $g$. Mumford conjectured that the stable cohomology
ring $H^*({\mathcal M}_\infty; \QQ)$ is freely generated by the tautological classes $\kappa_1, \kappa_2, \ldots$ and Miller
\cite{Miller} and Morita \cite{Morita1987} proved part of this conjecture by showing that the
natural map
$$ \QQ[\kappa_1, \kappa_2, \ldots ] \to 
H^*({\mathcal M}_\infty; \QQ) $$
is injective. The remainder of Mumford's conjecture, that this
map is surjective, remained
unproved for nearly two decades. However Madsen and Tillmann
\cite{MT} found an interpretation of Mumford's map 
on the level of homotopy, which they conjectured should be
a homotopy equivalence.  
Very recently 
a proof of their conjecture,
using h-principle arguments combined with 
Harer stabilisation, has been announced by
Madsen and Weiss \cite{MadsenWeiss,Tillmann}, and from this Mumford's conjecture follows.

The tautological ring of ${{\mathcal M}}_g$ for finite $g$ 
has many beautiful properties. 
Faber \cite{Faber} conjectured
that when $g \geq 2$ the tautological ring of 
${{\mathcal M}}_g$ looks like the algebraic
cohomology ring of a nonsingular complex projective variety
of dimension $g-2$, 
 and that it is generated by the tautological classes $\kappa_1, \kappa_2, \ldots,
\kappa_{[g/3]}$ with no relations in degrees at most $[g/3]$.
He also provided an explicit conjecture for a complete set of relations
among these generators. 
Progress has been made by many contributors towards
 Faber's conjectures, 
and also 
related problems on moduli spaces linked to ${{\mathcal M}}_g$.  
In particular Morita \cite{Morita, Morita2} has recently proved that the 
rational cohomological version of the tautological ring of
${{\mathcal M}}_g$ is indeed generated by 
$\kappa_1, \kappa_2, \ldots,
\kappa_{[g/3]}$. The definition of the tautological
ring has also been extended to 
the compactification $\overline{{\mathcal M}}_{g,n}$ of
the moduli space ${{\mathcal M}}_{g,n}$
of nonsingular curves of genus $g$ with $n$ marked points
(motivated by Witten's conjectures \cite{Witten},
proved by Kontsevich \cite{Kontsevich}, on intersection
pairings on $\overline{{\mathcal M}}_{g,n}$). 

The moduli spaces ${{\mathcal M}}_{g}$ and 
${{\mathcal M}}_{g,n}$ have other younger and more sophisticated relatives, such as the moduli spaces
${\mathcal M}_{g,n}(X,\beta)$ which parametrise holomomorphic
maps $f:\Sigma \to X$ from a nonsingular complex projective curve
$\Sigma$ of genus $g$ with $n$ marked points satisfying $f_*[\Sigma]=\beta
\in H_2(X)$, and their compactifications
$\overline{{\mathcal M}}_{g,n}(X,\beta)$ which parametrise \lq stable' maps. Intersection theory on $\overline{{\mathcal M}}_{g,n}(X,\beta)$
is fundamental to Gromov-Witten theory and quantum cohomology for $X$, 
with numerous applications in the last decade to enumerative
geometry. The Virasoro conjecture of Eguchi, Hori and Xiong
provides relations among the descendent Gromov-Witten invariants of $X$, 
and its recent proof by Givental \cite{Givental}
for $X=\PP^n$ implies part of Faber's conjecture by \cite{GP}.

Other relatives of ${\mathcal M}_g$ include the moduli spaces of pairs $(\Sigma, E)$ 
where $\Sigma$ is a nonsingular curve and $E$ is a stable vector
bundle over $\Sigma$, and their compactifications; intersection theory
on these relates intersection theory on $\overline{{\mathcal M}}_g$
and intersection theory on moduli spaces of bundles over a fixed curve,
which is by now quite well understood. 

\renorm
\section{Moduli spaces of curves}

The study of algebraic curves, and how they vary in families, has been fundamental
to algebraic geometry since the beginning of the subject, and has made
huge advances in the last few decades \cite{ACGH,HarrisMorrison}. 
The concept of moduli as parameters
describing as efficiently as possible the variation of geometric objects
was initiated in Riemann's famous paper \cite{Riemann} of 1857, in
which he observed that an isomorphism class of compact Riemann surfaces of
genus $g \geq 2$ \lq h\"{a}ngt ... von $3g-3$ stetig ver\"{a}nderlichen Gr\"{o}ssen
ab, welche die Moduln dieser Klasse genannt werden sollen'. In modern terminology,
Riemann's observation is the statement that the dimension of 
${\mathcal M}_g$ is $3g-3$ if $g \geq 2$. It was not until the 1960s that precise
definitions and methods of constructing moduli spaces were given by Mumford in
\cite{GIT} following ideas of Grothendieck. 
Roughly speaking, the moduli space ${\mathcal M}_g$ is the set of isomorphism classes of
nonsingular complex projective curves\footnote{All complex curves and
real surfaces will be assumed to be connected.} of genus $g$, endowed with the structure of
a complex variety in such a way that any family of nonsingular complex projective
curves parametrised by a base space $S$ induces a morphism from $S$ to ${\mathcal M}_g$
which associates to each $s \in S$ the isomorphism class of the curve parametrised
by $s$. The moduli spaces ${{\mathcal M}}_g$ can be constructed in several different ways,
including
\begin{itemize}
\item as orbit spaces for group actions,
\item via period maps and Torelli's theorem, and
\item using Teichm\"{u}ller theory.
\end{itemize}

The first of these is a standard method for constructing many different
moduli spaces, using Mumford's geometric invariant theory \cite{GIT,N2,Viehweg}
or more recent ideas due to Koll\'{a}r \cite{Kollar} and to Mori and Keel \cite{KeelMori}.
Geometric invariant theory provides a beautiful compactification of ${\mathcal M}_g$ known
as the Deligne-Mumford compactification $\overline{{\mathcal M}}_g$ \cite{DeligneMumford}.
This compactification is itself modular: it is the moduli space of (Deligne-Mumford) stable
curves (i.e. complex projective curves with only nodal singularities and
finitely many automorphisms). 
$\overline{{\mathcal M}}_{g}$ is singular but in a
relatively mild way; it is the quotient of a nonsingular variety by a finite group
action \cite{Looijenga1994}. 

The moduli space ${{\mathcal M}}_{g,n}$ of nonsingular
complex projective curves of genus $g$ with $n$ marked points has a similar
compactification $\overline{{\mathcal M}}_{g,n}$ which is the moduli space
of complex projective curves with $n$ marked nonsingular points and with only nodal
singularities and finitely many automorphisms. Finiteness of the automorphism
group of such a curve $\Sigma$ is equivalent to the requirement that any irreducible component of genus
0 (respectively 1) has at least 3 (respectively 1) special points, where \lq special'
means either marked or singular in $\Sigma$ (and the condition on genus 1 components here is
redundant when $g\geq 2$).

The second method of construction using the period
matrices of curves leads to a different compactification $\tilde{{\mathcal M}}_g$
of ${\mathcal M}_g$ known as the Satake (or Satake-Baily-Borel) compactification. Like the Deligne-Mumford
compactification, $\tilde{{\mathcal M}}_g$ is a complex projective variety, but the
boundary $\tilde{{\mathcal M}}_g \setminus {\mathcal M}_g$ of ${\mathcal M}_g$ in
$\tilde{{\mathcal M}}_g$ has (complex) codimension 2 for $g \geq 3$ whereas the 
boundary $\Delta = \overline{{\mathcal M}}_g \setminus {\mathcal M}_g$ of 
${{\mathcal M}}_g$ in
$\overline{{\mathcal M}}_g$ has codimension 1. Each of the irreducible 
components $\Delta_0, \ldots,
\Delta_{[g/2]}$ of $\Delta$  is the closure of a locus of curves with exactly one node
(irreducible curves with one node in the case of $\Delta_0$, and in the
case of any other $\Delta_i$ the union
of two nonsingular curves of genus $i$ and $g-i$ meeting at a single point).
The divisors $\Delta_i$ meet transversely in $\overline{{\mathcal M}}_g$,
and their intersections define a natural decomposition of $\Delta$ into
connected strata which parametrise stable curves of a fixed topological
type. 
The boundary of ${{\mathcal M}}_{g,n}$ in $\overline{{\mathcal M}}_{g,n}$
has a similar description, but now as well as the genus of each irreducible
component it is necessary to keep track of which marked
points it contains. 

The third method of constructing ${{\mathcal M}}_g$, via Teichm\"{u}ller theory,
leaves algebraic geometry altogether.  

\renorm
\section{Teichm\"{u}ller theory and mapping class groups}

Important recent advances concerning the cohomology of the moduli spaces ${{\mathcal M}}_g$
(in particular \cite{MT, MadsenWeiss, Morita, Morita2,Tillmann}) have been proved by
topologists via the link between these moduli spaces and mapping class groups
of compact surfaces.

Let us fix a compact oriented smooth surface $\Sigma_g$ of genus
$g \geq 2$, and let $\mbox{Diff}_+ \Sigma_g$ be the group of orientation preserving
diffeomorphisms of $\Sigma_g$. Then the mapping class group $\Gamma_g$ of $\Sigma_g$
is the group
$$ \Gamma_g = \pi_0 (\mbox{Diff}_+ \Sigma_g)$$
of connected components of $\mbox{Diff}_+ \Sigma_g$. It acts properly and discontinuously
on the Teichm\"{u}ller space ${{\mathcal T}}_g$ of $\Sigma_g$, which is the space of conformal
structures on $\Sigma_g$ up to isotopy. The Teichm\"{u}ller space
${{\mathcal T}}_g$ is homeomorphic to $\RR^{6g-6}$, and its quotient by the action of
the mapping class group $\Gamma_g$ can be identified naturally with the moduli
space ${{\mathcal M}}_g$. This means that there is a natural isomorphism of
rational cohomology 
\begin{equation} \label{1} H^*({{\mathcal M}}_g;\QQ) \cong H^*(\Gamma_g;\QQ).\end{equation}
The corresponding integral cohomology groups are not in general isomorphic because of the
existence of nonsingular complex projective curves with nontrivial automorphisms.
If, however, we work with the moduli spaces ${{\mathcal M}}_{g,n}$ of nonsingular
complex projective curves of genus $g$ with $n$ marked points, then when $n$ is
large enough such marked curves have no nontrivial automorphisms (cf. \cite{HarrisMorrison}
p 37) and
$$H^*({{\mathcal M}}_{g,n};\ZZ) \cong H^*(\Gamma_{g,n};\ZZ)$$
where $\Gamma_{g,n}$ is the group of connected
components of the group $\mbox{Diff}_+ \Sigma_{g,n}$ of orientation preserving
diffeomorphisms of $\Sigma_g$ which fix $n$ chosen points on $\Sigma_g$.

In fact \cite{EarleEels} the components of $\mbox{Diff}_+ \Sigma_g$ are contractible when
$g \geq 2$, so there is also a natural isomorphism 
$$H^*(\Gamma_g;\ZZ) \cong H^*(B \mbox{Diff}_+ \Sigma_g;\ZZ)$$
where $B\mbox{Diff}_+ \Sigma_g$ is the universal classifying space for
$\mbox{Diff}_+ \Sigma_g$. This means that any cohomology class of the mapping class
group $\Gamma_g$ can be regarded as a characteristic class of oriented surface bundles,
while any rational cohomology class of $\Gamma_g$ can be regarded as a rational cohomology
class of the moduli space ${{\mathcal M}}_g$.

The mapping class group $\Gamma_g$ can be described in a group theoretical
way. $\Gamma_g$ acts faithfully by outer automorphisms (that is, the action
is defined modulo inner automorphisms) on the fundamental group $\pi_1(\Sigma_g)$
of $\Sigma_g$,
which is generated by $2g$ elements $a_1,\ldots,a_g,b_1,\ldots,b_g$ subject to one
relation $\prod_j a_j b_j a_j^{-1} b_j^{-1} =1$, and the image of $\Gamma_g$ in
$\mbox{Out}(\pi_1(\Sigma_g))$ is the group of outer automorphisms of $\pi_1(\Sigma_g)$
which act trivially on $H_2(\pi_1(\Sigma_g);\ZZ)$ \cite{ZVC}. $\Gamma_g$ 
has a finite presentation \cite{Wajnryb} with generators represented by Dehn twists
(diffeomorphisms of $\Sigma_g$ obtained by cutting $\Sigma_g$ along a regularly
embedded circle, twisting one of the resulting boundary circles through $2\pi$ and
reglueing). There are similar descriptions of $\Gamma_{g,n}$ \cite{ZVC,Gervais}.

\renorm
\section{Stable cohomology}

Harer \cite{Harer} proved in the 1980s that $H^k(\Gamma_g;\ZZ)$ and $H^k(\Gamma_{g+1};\ZZ)$
are isomorphic when $g \geq 3k-1$, and the same is true for $\Gamma_{g,n}$ and $\Gamma_{g+1,n}$.
This bound was improved by Ivanov \cite{Ivanov} and Harer \cite{Harerpre}
made a further improvement. Since $H^*(\Gamma_g;\QQ)$ is isomorphic to $H^*({{\mathcal M}}_g;
\QQ)$, this means that the rational cohomology group $H^k ({{\mathcal M}}_g;\QQ)$ is independent of $g$ for $g \gg k$, and we can define the stable cohomology ring 
$$H^*({{\mathcal M}}_\infty ; \QQ)$$
so that $H^k({{\mathcal M}}_\infty;\QQ) \cong H^k({{\mathcal M}}_g;\QQ)$ for $g \gg k$.

Harer's stabilisation map can be defined as follows. We choose a smooth identification of
$\Sigma_{g+1}$ with a connected sum 
of a smooth surface $\Sigma_g$
of genus $g$ and a surface $\Sigma_1$ of genus 1 (and if we have marked points we make
sure they all correspond to points in $\Sigma_g$). Let $\Gamma_{g+1,\Sigma_1}$ be the
subgroup of $\Gamma_{g+1}$ consisting of mapping classes represented by diffeomorphisms
from $\Sigma_{g+1}$ to itself which fix all the points coming from $\Sigma_1$. The
result of collapsing all such points in $\Sigma_{g+1}$ together is diffeomorphic to
$\Sigma_g$, so there is a homomorphism from $\Gamma_{g+1,\Sigma_1}$ to $\Gamma_g$ as
well as an inclusion of $\Gamma_{g+1,\Sigma_1}$ in $\Gamma_{g+1}$. Harer showed that
both of these induce isomorphisms 
$$H^k(\Gamma_g;\ZZ) \cong H^k(\Gamma_{g+1,\Sigma_1};\ZZ) \mbox{ and }
H^k(\Gamma_{g+1};\ZZ) \cong H^k(\Gamma_{g+1,\Sigma_1};\ZZ)$$
when $g \gg k$, and likewise we have 
\begin{equation} \label{11} H^k(\Gamma_{g,n};\ZZ) \cong H^k(\Gamma_{g+1,\Sigma_1,n};\ZZ) \mbox{ and }
H^k(\Gamma_{g+1,n};\ZZ) \cong H^k(\Gamma_{g+1,\Sigma_1,n};\ZZ) \end{equation}
when $g \gg k$. 

A similar construction can be made to describe the stabilisation isomorphism
$$ H^k({{\mathcal M}}_{g,n};\QQ) \cong H^k({{\mathcal M}}_{g+1,n};\QQ)$$
for the moduli spaces ${{\mathcal M}}_{g,n}$ (cf. [28] p31). Identifying the last marked
point of a smooth nonsingular complex projective curve of genus $g$ with a marked
point on a curve of genus 1 gives a stable curve of genus $g+1$ with $n$ marked
points. This defines for us a morphism 
$$\phi: {{\mathcal M}}_{g,n+1} \times {{\mathcal M}}_{1,1} \to 
\overline{{\mathcal M}}_{g+1,n}$$
whose image is an open subset of an irreducible component of the boundary
of ${{\mathcal M}}_{g+1,n}$ in $\overline{{\mathcal M}}_{g+1,n}$, and
there is a normal bundle ${\mathcal N}_\phi$ which is a complex line bundle
(in the sense of orbifolds) over ${{\mathcal M}}_{g,n+1} \times {{\mathcal M}}_{1,1}$.
Using ${\mathcal N}_\phi^*$ to denote the complement of the zero section of
${\mathcal N}_\phi$ we can compose projection maps with the forgetful map
from ${{\mathcal M}}_{g,n+1}$ to ${{\mathcal M}}_{g,n}$ to get
$${\mathcal N}^*_\phi \to {{\mathcal M}}_{g,n+1} \times {{\mathcal M}}_{1,1}
\to {{\mathcal M}}_{g,n+1} \to {{\mathcal M}}_{g,n}$$
which induces
\begin{equation} \label{2} H^k({{\mathcal M}}_{g,n};\QQ) \to H^k({\mathcal N}_\phi^*;
\QQ).\end{equation}
On the other hand, using a tubular neighbourhood of the image of $\phi$ in
$\overline{{\mathcal M}}_{g+1,n}$ we obtain a natural homotopy class of maps from
${\mathcal N}^*_\phi$ to ${{\mathcal M}}_{g+1,n}$ which induces
\begin{equation} \label{3} H^k({{\mathcal M}}_{g+1,n};\QQ) \to H^k({\mathcal N}_\phi^*;
\QQ).\end{equation}
Here (\ref{2}) and (\ref{3}) represent Harer's maps (\ref{11}) and hence they are isomorphisms
if $g \gg k$.

\renorm
\section{Tautological classes}

When $g\geq 2$ Mumford \cite{Mumford} and Morita \cite{Morita1984} independently defined tautological
classes 
$$\kappa_i \in H^{2i}(\overline{{\mathcal M}}_{g};\QQ) \mbox{ and }e_i \in H^{2i}(\Gamma_g;\ZZ)$$ 
which correspond up to a sign $(-1)^{i+1}$ in $H^*({{\mathcal M}}_{g};\QQ)$ 
under the isomorphism (\ref{1}). The subalgebra $R^*({\mathcal M}_g)$ 
of $H^*({{\mathcal M}}_{g};\QQ)$
generated by the $\kappa_i$, or equivalently by the $e_i$, is called its
tautological ring.

The classes $\kappa_i$ are defined using the natural
forgetful map $\pi:{{\mathcal M}}_{g,1} \to {{\mathcal M}}_{g}$ which takes an element
$[\Sigma, p]$
of ${{\mathcal M}}_{g,1}$ represented by a nonsingular complex projective curve 
$\Sigma$ with one marked point $p$ to the element $[\Sigma]$ of ${{\mathcal M}}_{g}$ represented
by $\Sigma$. This is often called the universal curve over ${{\mathcal M}}_{g}$, since
for generic choices of $\Sigma$ the fibre $\pi^{-1}([\Sigma])$ is a copy of $\Sigma$.
However if $\Sigma$ has nontrivial automorphisms then $\pi^{-1}([\Sigma])$ is not
a copy of $\Sigma$ but is instead the quotient of $\Sigma$ by its automorphism group
$\mbox{Aut}(\Sigma)$ (which has size at most $84(g-1)$ when $g \geq 2$).

From the topologists' viewpoint the r\^{o}le of $\pi:{\mathcal M}_{g,1}
\to {\mathcal M}_g$ is played by the universal oriented $\Sigma_g$-bundle
$$\Pi : E \mbox{Diff}_+ \Sigma_g \to B \mbox{Diff}_+ \Sigma_g.$$
Its relative tangent bundle is an oriented real vector bundle of rank 2 on
$E \mbox{Diff}_+ \Sigma_g$ (whose fibre at $x \in E \mbox{Diff}_+ \Sigma_g$
is the tangent space at $x$ to the oriented surface $\Pi^{-1}(x)$), so it
has an Euler class $e \in H^2(E\mbox{Diff}_+ \Sigma_g;\ZZ)$. Morita defined
his tautological classes
$$e_i \in H^{2i}(\Gamma_g;\ZZ) \cong H^{2i}( B \mbox{Diff}_+ \Sigma_g ;\ZZ)$$
by setting $e_i$ to be the pushforward (or integral over the fibres) 
$\Pi_!(e^{i+1})$ of $e^{i+1}$.

To define his tautological classes $\kappa_i$ Mumford used essentially the
same procedure with the forgetful map $\pi:{{\mathcal M}}_{g,1} \to {{\mathcal M}}_{g}$,
except that he used cotangent spaces instead of tangent spaces (which is the reason
that $\kappa_i$ and $e_i$ only correspond up to a sign $(-1)^{i+1}$) and the relative
cotangent bundle (or relative dualising sheaf) for $\pi:{{\mathcal M}}_{g,1} \to {{\mathcal M}}_{g}$ exists as a complex line bundle over ${{\mathcal M}}_{g,1}$ only in the
sense of orbifold line bundles (or line bundles over stacks) because of the
existence of nontrivial automorphism groups $\mbox{Aut}(\Sigma)$. 

The forgetful map $\pi:{{\mathcal M}}_{g,1} \to {{\mathcal M}}_{g}$
can be generalised to $\pi:{{\mathcal M}}_{g,n+1} \to {{\mathcal M}}_{g,n}$
for any $n \geq 0$ by forgetting the last marked point of an $n+1$-pointed curve, and
this can be extended to $\pi:\overline{{\mathcal M}}_{g,n+1} \to 
\overline{{\mathcal M}}_{g,n}$. Care is needed here when the last marked point
lies on an irreducible component with genus 0 and only two other special
points; such an irreducible component needs to be collapsed in order to produce
a stable $n$-pointed curve of genus $g$. This collapsing procedure gives us
a forgetful map $\pi:\overline{{\mathcal M}}_{g,n+1} \to 
\overline{{\mathcal M}}_{g,n}$ whose fibre at $[\Sigma, p_1, \ldots, p_n] \in
\overline{{\mathcal M}}_{g,n}$ can be identified with the quotient of $\Sigma$
by the automorphism group of $(\Sigma, p_1,...,p_n)$. Mumford's tautological classes can be extended to classes
$\kappa_i \in H^{2i}(\overline{{\mathcal M}}_{g,n};\QQ)$ (in fact
to classes in the rational Chow ring of $\overline{{\mathcal M}}_{g,n}$) defined
by
$$\kappa_i = \pi_! (c_1(\omega_{g,n})^{i+1})$$
where $\omega_{g,n}$ is the relative dualising sheaf of  
$\pi:\overline{{\mathcal M}}_{g,n+1} \to 
\overline{{\mathcal M}}_{g,n}$ and $c_1(\omega_{g,n})\in H^2(
\overline{{\mathcal M}}_{g,n};\QQ)$ is its first Chern class.

When $n>0$ there are other interesting tautological classes on
${{\mathcal M}}_{g,n}$ and $\overline{{\mathcal M}}_{g,n}$
exploited by Witten. 
The forgetful map $\pi:\overline{{\mathcal M}}_{g,n+1} \to 
\overline{{\mathcal M}}_{g,n}$ has tautological sections $s_j:
\overline{{\mathcal M}}_{g,n} \to 
\overline{{\mathcal M}}_{g,n+1}$ for $1 \leq j \leq n$ such that
$s_j([\Sigma,p_1, \ldots, p_n])$ is the element of 
$\pi^{-1}([\Sigma,p_1, \ldots, p_n]) = \Sigma/\mbox{Aut}(\Sigma)$ represented
by $p_j$.
The Witten classes $\psi_j \in H^2(\overline{{\mathcal M}}_{g,n};\QQ)$
for $j=1, \ldots,n$ can then be defined by
$$\psi_j = c_1(s^*_j(\omega_{g,n})).$$
Roughly speaking, $\psi_j$ is the first Chern class of the (orbifold)
line bundle on $\overline{{\mathcal M}}_{g,n}$ whose fibre at 
$[\Sigma, p_1, \ldots, p_n]$ is the cotangent space $T^*_{p_j}\Sigma$
to $\Sigma$ at $p_j$. 

The boundary $\Delta = \overline{{\mathcal M}}_{g,n} \setminus {{\mathcal M}}_{g,n}$
of ${{\mathcal M}}_{g,n}$ in $\overline{{\mathcal M}}_{g,n}$ is the union of
finitely many divisors which meet transversely in $\overline{{\mathcal M}}_{g,n}$.
The intersection of any nonempty set of these divisors is the closure of
a subset of ${{\mathcal M}}_{g,n}$ parametrising stable $n$-pointed curves of
some fixed topological type, and is the image of a finite-to-one map to
$\overline{{\mathcal M}}_{g,n}$ from a product of moduli spaces
$\prod_k \overline{{\mathcal M}}_{g_k,n_k}$ which glues together stable
curves of genus $g_k$ with $n_k$ marked points at certain of the
marked points. These glueing maps induce pushforward maps on cohomology
\begin{equation} \label{4} H^*(\prod_k \overline{{\mathcal M}}_{g_k,n_k};\QQ) \to
H^*(\overline{{\mathcal M}}_{g,n};\QQ) \end{equation}
and the tautological ring $R^*(\overline{{\mathcal M}}_{g,n};\QQ)$
is defined inductively to be the subalgebra of $H^*(\overline{{\mathcal M}}_{g,n};\QQ)$
generated by the Mumford classes, the Witten classes and the
images of the tautological classes in 
$H^*(\prod_k \overline{{\mathcal M}}_{g_k,n_k};\QQ)$ under the pushforward
maps (\ref{4}) from the boundary of ${{\mathcal M}}_{g,n}$. Its restriction
to $H^*({{\mathcal M}}_{g,n};\QQ)$ is the tautological ring of
${{\mathcal M}}_{g,n}$ and is generated by the Mumford and Witten classes.

\renorm
\section{Mumford's conjecture}

Mumford's tautological classes $\kappa_i \in H^{2i}({{\mathcal M}}_{g};\QQ)$
are preserved by Harer stabilisation when $g$ is sufficiently large, and so
they define elements of the stable cohomology $H^*({{\mathcal M}}_{\infty};\QQ)$.
Mumford conjectured in \cite{Mumford} that $H^*({{\mathcal M}}_{\infty};\QQ)$ is
freely generated by $\kappa_1, \kappa_2, \ldots$, or in other words that the obvious map
\begin{equation} \label{mumf} \QQ[\kappa_1,\kappa_2,\ldots ] \to H^*({{\mathcal M}}_{\infty}
;\QQ) \end{equation}
is an isomorphism. Miller \cite{Miller} and Morita \cite{Morita1987} soon proved that
this map is injective, so it remained to prove surjectivity. Not long ago Madsen and 
Tillmann \cite{MT} found a homotopy version of Mumford's map (\ref{mumf}) which they
conjectured to be a homotopy equivalence, and very recently
Madsen and Weiss \cite{MadsenWeiss} have announced
a proof of their conjecture, from which Mumford's conjecture follows.

The Madsen-Tillmann map involves the stable mapping class group $\Gamma_\infty$ rather
than the moduli spaces ${{\mathcal M}}_{g}$. From the description of ${{\mathcal M}}_{g}$
as the quotient of the $\Gamma_g$ action on Teichm\"{u}ller space ${\mathcal T}_g$ it
follows that when $g \geq 2$ there is a continuous map
\begin{equation} \label{7} B \Gamma_g \to {{\mathcal M}}_{g} \end{equation}
uniquely determined up to homotopy. It is known that $\Gamma_g$ is a perfect group
when $g \geq 3$ \cite{Harer1983}, so we can apply Quillen's plus construction
to $B\Gamma_g$ to obtain a simply connected space $B \Gamma_g^+$ with the same homology as
$B \Gamma_g$. The moduli space ${{\mathcal M}}_{g}$ is also simply connected, so (\ref{7})
factors through a map 
$B \Gamma_g^+ \to {{\mathcal M}}_{g}$
which induces the isomorphism
$H^*(\Gamma_g;\QQ) \to H^*({{\mathcal M}}_{g};\QQ)$
 discussed above at (\ref{1}). Moreover Harer stabilisation gives us maps
$B\Gamma_g^+ \to B\Gamma_{g+1}^+$ between simply connected spaces which are
homology equivalences (and hence also homotopy equivalences) in a range up to
some degree which tends to infinity with $g$. If $B\Gamma_\infty^+$ denotes the
homotopy direct limit of these maps as $g \to \infty$, then Mumford's conjecture
becomes the statement that 
$$H^*(B\Gamma_\infty^+;\QQ) \cong \QQ[\kappa_1, \kappa_2, \ldots ].$$
The conjecture of Madsen and Tillmann \cite{MT} describes the homotopy type of
$B \Gamma_\infty^+$ (or rather $\ZZ \times B \Gamma_{\infty}^+$), 
giving Mumford's conjecture as a corollary.

Tillmann \cite{Tillmann} had already shown that $\ZZ \times B \Gamma_\infty^+$
is an infinite loop space, in the sense that there exists a sequence of spaces
$E_n$ with $E_n = \Omega E_{n+1}$ and $\ZZ \times B\Gamma_\infty^+ = E_0$.
This was an encouraging result because infinite loop spaces have many good
properties. Subsequently Madsen and Tillmann \cite{MT} found an $\Omega^\infty$ map 
$\alpha_\infty$ from
$\ZZ \times B \Gamma_\infty^+$ to an infinite loop space which they denoted by
$\Omega^\infty \CC \PP^\infty_{-1}$ and whose connected
component has rational
cohomology isomorphic to $\QQ[\kappa_1,\kappa_2,\ldots]$.

The infinite loop space $\Omega^\infty \CC \PP^\infty_{-1}$ is related
to the limit $\CC \PP^\infty$ of the complex projective spaces $\CC\PP^k$
as $k \to \infty$. Over $\CC\PP^k$ there is a tautological
complex line bundle $L_k$, whose fibre at $x \in \CC\PP^k$ is the one-dimensional
subspace of $\CC^{k+1}$ represented by $x$, and a complex vector bundle
$L_k^\bot$ of rank $k$ which is its complement in the trivial bundle of
rank $k+1$ over $\CC\PP^k$. The restriction of $L_{k+1}^\bot$ to $\CC\PP^k$
is the direct sum of $L_k^\bot$ and a trivial complex line bundle, giving
us maps 
$\mbox{Th}(L_k^\bot) \to \Omega^2 \mbox{Th}(L^\bot_{k+1})$ and 
$\Omega^{2k+2}\mbox{Th}(L_k^\bot) \to \Omega^{2k+4} \mbox{Th}(L^\bot_{k+1})$
where $\mbox{Th}(L_k^\bot)$ is the Thom space (or one-point compactification)
of the bundle $L_k^\bot$. Madsen and Tillmann define $\Omega^\infty \CC \PP^\infty_{-1}$
to be the direct limit of the spaces $\Omega^{2k+2}\mbox{Th}(L_k^\bot)$ as 
$k \to \infty$.

Homotopy classes of maps from an $n$-dimensional manifold $X$ to $\Omega^\infty
\CC\PP^\infty_{-1}$ are represented by proper maps $\phi:M \to X$ from an
$(n+2)$-dimensional manifold $M$ together with an \lq artificial differential'
$\Phi:TM \to \phi^* TX$ and an orientation of $\mbox{ker}\Phi$. Here
$\Phi$ is a stable vector bundle surjection; that is, it may be that $\Phi$
is defined and becomes a surjective bundle map only once a
trivial bundle of sufficiently large rank has been added to $TM$ and $\phi^*TX$.
Any smooth oriented surface
bundle $\phi:E \to X$ induces a homotopy class of maps from $X$ to $\Omega^\infty
\CC\PP^\infty_{-1}$ represented by $\phi$ together with its differential
$\Phi = d\phi:TE \to \phi^* TX$, and this effectively defines the Madsen-Tillmann
map $\alpha_\infty: \ZZ \times B \Gamma^+_\infty \to \Omega^\infty
\CC\PP^\infty_{-1}$.

Submersion theory suggests a way to tackle the problem of showing that $\alpha_\infty$
is a homotopy equivalence, but compactness of $X$ creates a difficulty for this.
Therefore Madsen and Weiss replace $X$ with $X \times \RR$. They
study a commutative diagram
$$\begin{array}{ccccc} {\mathcal V} & \rightarrow & {\mathcal W} &
\rightarrow & {\mathcal W}_{\mbox{loc}} \\
\downarrow & & \downarrow & & \downarrow \\
h{\mathcal V} & \rightarrow & h{\mathcal W} & \rightarrow & h{\mathcal W}_{\mbox{loc}}
\end{array}$$
of contravariant functors from smooth manifolds to sets with the sheaf property for
open coverings, and the induced diagram
$$\begin{array}{ccccc} |{\mathcal V}| & \rightarrow & |{\mathcal W}| &
\rightarrow & |{\mathcal W}_{\mbox{loc}}| \\
\downarrow & & \downarrow & & \downarrow \\
|h{\mathcal V}| & \rightarrow & |h{\mathcal W}| & \rightarrow & |h{\mathcal W}_{\mbox{loc}}|
\end{array}$$
of the associated spaces, where homotopy classes of maps from
$X$ to $|{\mathcal F}|$ correspond naturally to concordance classes in ${\mathcal F}(X)$,
and $s_0,s_1 \in {\mathcal F}(X)$ are concordant if $s_0 = t|_{X \times \{ 0\}}$
and $s_1 = t|_{X \times \{1\}}$ for some $t \in {\mathcal F}(X \times \RR)$.

If $X$ is any smooth manifold then elements of ${\mathcal V}(X)$
are given by smooth oriented surface bundles $E$ 
(that is, proper submersions whose fibres are {connected}
oriented surfaces) over $X\times \RR$, together with identifications 
$\partial E \cong \partial (S^1 \times [0,1] \times X \times \RR)$ compatible
with the maps to $X \times\RR$. These identifications on the boundary are
crucial, because they give ${\mathcal V}$ and the other functors involved the
structure of monoids, and thus the associated spaces become topological
monoids. 

In one version of the bottom row $
h{\mathcal V}  \rightarrow  h{\mathcal W}  \rightarrow  h{\mathcal W}_{\mbox{loc}}$
of the commutative diagram, 
elements of $h{\mathcal V}(X)$
are given by $(n+3)$-dimensional manifolds $E$, where $n=\mbox{dim}X$, and smooth maps $\pi:E \to X$ and
$f,g:E \to \RR$ such that $(\pi,f):E \to X \times \RR$ is a submersion and
$(\pi,g):E \to X \times \RR $ is proper, together with an identification
$\partial E \cong \partial (S^1 \times [0,1] \times X \times \RR)$ compatible
with the maps to $X$ and $\RR$. 
If $(\pi,f):E \to X \times \RR$ represents an element of ${\mathcal V}(X)$
then we get an element of $h{\mathcal V}(X)$ by setting $g=f$.
The functors 
${\mathcal W}$ and $h{\mathcal W}$ are defined similarly, except that the requirement
that $(\pi,f):E \to X \times \RR$ should be a submersion is weakened to the
requirements that $\pi:E \to X$ should be a submersion and that the restriction of
$f:E \to \RR$ to any fibre of $\pi$ should be a Morse function. For
${\mathcal W}_{\mbox{loc}}$ and $h{\mathcal W}_{\mbox{loc}}$ the
requirements are weakened again,
so that \lq proper' is replaced by \lq proper when restricted to the set of
singularities of $f$ on fibres of $\pi$'.

The strategy of Madsen and Weiss is to deduce that $\alpha_\infty$ is a homotopy
equivalence from the following properties of the commutative diagram above:

\noindent (i) the first vertical map represents the Madsen-Tillmann map $\alpha_\infty$;

\noindent (ii) the second vertical map is a homotopy equivalence (by a corollary to
 Vassiliev's
h-principle \cite{Vassiliev});

\noindent (iii) the third vertical map is also a homotopy equivalence (by a much easier
argument);

\noindent (iv) the bottom row is a homotopy fibre sequence;

\noindent (v) the top row becomes a homotopy fibre sequence after group completion
(using stratifications of $|{\mathcal W}|$ and
$|{\mathcal W}_{\mbox{loc}}|$ and a subtle application of Harer stabilisation).

\renorm
\section{Faber's conjectures}

Although Mumford's conjecture tells us that the tautological classes $\kappa_i$ generate
the stable cohomology ring $H^*({\mathcal M}_\infty;\QQ)$, they do not generate 
$H^*({\mathcal M}_g;\QQ)$ for finite $g$, and in fact $H^*({\mathcal M}_g;\QQ)$ has
lots of unstable cohomology (at least when $g$ is large enough). This follows from the
calculation  of 
Euler characteristics by Harer and Zagier \cite{HarerZagier} (see also \cite{Kontsevich}).
They show that the orbifold Euler characteristic of ${{\mathcal
M}}_{g,n}$ is 
$$(-1)^{n-1} \frac{(2g+n-3)!}{(2g-2)!} \zeta(1-2g)$$
where $\zeta$ denotes the Riemann $\zeta$-function, and their work implies that when $g \geq 15$
the Euler characteristic of ${{\mathcal M}}_g$ is too large in absolute value for
$H^*({{\mathcal M}}_g;\QQ)$ to be generated by $\kappa_1,\kappa_2,\ldots$ (cf. also
\cite{GraberPandharipande2001,Looijenga1993}). Nonetheless the tautological ring $R^*({{\mathcal M}}_g)$ generated by $\kappa_1,\kappa_2,\ldots$ has many beautiful properties.

Faber \cite{Faber} has conjectured that $R^*({{\mathcal M}}_g)$ has the structure of the
algebraic cohomology ring of a nonsingular complex projective variety of dimension
$g-2$. More precisely, he conjectured that

\noindent (i) $R^k({\mathcal M}_g)$ is zero when $k>g-2$ and is one-dimensional when
$k=g-2$, and the natural pairing $R^k({\mathcal M}_g) \times R^{g-2-k}({\mathcal M}_g)
\to R^{g-2}({\mathcal M}_g)$ is perfect. In addition $R^k({\mathcal M}_g)$ satisfies
the Hard Lefschetz property and the Hodge index theorem with respect to the class 
$\kappa_1$.

\noindent(ii) The classes $\kappa_1, \ldots , \kappa_{[g/3]}$ generate $R^*({\mathcal M}_g)$
with no relations in degrees up to and including $[g/3]$.

\noindent(iii) Faber also gave an explicit conjecture for a complete set of relations
between these generators (in terms of the proportionalities between monomials in
$R^{g-2}({\mathcal M}_g)$).

When $g \leq 15$ Faber \cite{Faber} has proved all these conjectures concerning $R^*({\mathcal M}_g)$, and for general $g$ Looijenga \cite{Looijenga} and Faber \cite{Fabernon}
have shown that $R^k({\mathcal M}_g)$ is zero when $k>g-2$ and is one-dimensional when
$k=g-2$. Their proofs apply to both the cohomological version and the
Chow ring version of $R^*({\mathcal M}_g)$. Using topological methods, Morita
\cite{Morita, Morita2} has recently proved that the classes $\kappa_1, \ldots,
\kappa_{[g/3]}$ generate the cohomological version of $R^*({\mathcal M}_g)$
(and the rest of (ii) then follows essentially from \cite{Harerpre}).

The mapping class group $\Gamma_g$ acts naturally on $H_1(\Sigma_g;\ZZ)$ in a way 
which preserves the intersection pairing. This representation gives us an
exact sequence of groups
$$1 \to {\mathcal I}_g \to \Gamma_g \to Sp(2g;\ZZ) \to 1$$
where ${\mathcal I}_g$ denotes the subgroup of $\Gamma_g$ which acts trivially
on $H_1(\Sigma_g;\ZZ)$ and which is called the Torelli group. In \cite{Johnson0,
Johnson1,Johnson2,Johnson3,Johnson4}
Johnson showed that ${\mathcal I}_g$ is finitely generated for $g\geq 3$ (in
contrast with the case $g=2$ \cite{Mess}), introduced a surjective homomorphism 
$$\tau: {\mathcal I}_g \to {\wedge^3 H_1(\Sigma_g;\ZZ)}/{H_1(\Sigma_g;\ZZ)}$$
whose kernel is the subgroup of $\Gamma_g$ generated by all Dehn twists along
separating embedded circles, and used $\tau$ to determine the abelianisation of
${\mathcal I}_g$. Morita \cite{Morita1993} extended the Johnson homomorphism $\tau$
to a representation
$$\rho_1: \Gamma_g \to (\mbox{$\frac{1}{2}$} {\wedge^3 H_1(\Sigma_g;\ZZ)}/{H_1(\Sigma_g;\ZZ)})
\rtimes Sp(2g;\ZZ)$$
of the mapping class group $\Gamma_g$. Via the cohomology of semi-direct products
this induces
$$\rho_1^*: \mbox{Hom}(\wedge^*U,\QQ)^{Sp(2g;\ZZ)} \to H^*(\Gamma_g;\QQ) \cong 
H^*({\mathcal M}_g;\QQ)$$
where 
$U = {\wedge^3 H_1(\Sigma_g;\QQ)}/{H_1(\Sigma_g;\QQ)},$
and the image of $\rho_1^*$ is the tautological ring $R^*({\mathcal M}_g)$
\cite{KM,Looijenga1996}. By finding suitable relations in 
$ \mbox{Hom}(\wedge^*U,\QQ)^{Sp(2g;\ZZ)}$ and exploiting the 
map $H_1(\Sigma_g;\QQ) \to H_1(\Sigma_{g-1};\QQ)$ induced by collapsing
a handle of $\Sigma_g$, Morita \cite{Morita, Morita2} is able to prove
that the classes $\kappa_1, \ldots,
\kappa_{[g/3]}$ generate the cohomological version of $R^*({\mathcal M}_g)$.

Faber, Getzler, Hain, Looijenga, Pandharipande, Vakil and others (cf. \cite{Faber1990,Faber1990b,Faber1997,FaberPandharipande,GraberVakil,
HainLooijenga,Looijenga}) have also made
conjectures about the structure of the tautological rings of the compact moduli
spaces $\overline{{\mathcal M}}_{g,n}$, which are generated not just by the
Mumford classes $\kappa_i$ but also by the Witten classes $\psi_j$ and the 
pushforwards of tautological classes from the boundary of $\overline{{\mathcal M}}_{g,n}$.
For example, it is expected that $R^*(\overline{{\mathcal M}}_{g,n})$ looks like
the algebraic cohomology ring of a nonsingular complex projective variety of
dimension $3g-3+n$, while Getzler has conjectured that if $g>0$ then the monomials
of degree $g$ or higher in the Witten classes $\psi_j$ should all come from
the boundary of $\overline{{\mathcal M}}_{g,n}$ (a cohomological version of this
has been proved by Ionel \cite{Ionel}), and Vakil has made a
closely related conjecture that any tautological class in $R^k(\overline{{\mathcal 
M}}_{g,n})$ with $k\geq g$ should come from classes 
supported on boundary strata corresponding to
stable curves with at least $k-g+1$ components of genus 0.

\renorm
\section{The Virasoro conjecture}

The geometry of a nonsingular complex projective variety $X$ can be studied by examining
curves in $X$. Intersection theory on moduli spaces of curves in $X$, or more
precisely moduli spaces of maps from curves to $X$, leads to Gromov-Witten theory and the quantum cohomology of $X$,
with numerous applications in the last decade to enumerative geometry (cf. \cite{CK,
FultonPandharipande,Kontsevich,
Kontsevich1995,KontsevichManin}).

Let us assume for simplicity that $2g-2+n>0$. For any $\beta \in H_2(X;\ZZ)$ there is a moduli space ${{\mathcal M}}_{g,n}(X,\beta)$
of $n$-pointed nonsingular complex projective curves $\Sigma$ of genus $g$
equipped with maps $f:\Sigma \to X$ satisfying $f_*[\Sigma]=\beta$. This moduli space
has a compactification $\overline{{\mathcal M}}_{g,n}(X,\beta)$ which classifies
\lq stable maps' of type $\beta$ from $n$-pointed curves of genus $g$ into $X$
\cite{FultonPandharipande}. Here a map $f:\Sigma \to X$ from an $n$-pointed complex
projective curve $\Sigma$ satisfying $f_*[\Sigma] = \beta$
is called stable if $\Sigma$ has only nodal singularities and $f:\Sigma \to X$ has
only finitely many automorphisms, or equivalently every irreducible component of $\Sigma$
of genus 0 (respectively genus 1) which is mapped to a single point in $X$ by $f$ contains
at least 3 (respectively 1) special points. The forgetful map from
${{\mathcal M}}_{g,n}(X,\beta)$ to ${{\mathcal M}}_{g,n}$ which sends $[\Sigma,p_1,
\ldots, p_n, \, f:\Sigma \to X]$ to $[\Sigma, p_1,\ldots, p_n]$ extends to a forgetful map
$\pi:\overline{{\mathcal M}}_{g,n}(X,\beta) \to \overline{{\mathcal M}}_{g,n}$ which
collapses components of $\Sigma$ with genus 0 and at most two special points.

Of course, when $X$ is itself a single point, ${{\mathcal M}}_{g,n}(X,\beta)$
and $\overline{{\mathcal M}}_{g,n}(X,\beta)$ are simply the moduli spaces
${{\mathcal M}}_{g,n}$ and $\overline{{\mathcal M}}_{g,n}$. In general
$\overline{{\mathcal M}}_{g,n}(X,\beta)$ has more serious singularities
than $\overline{{\mathcal M}}_{g,n}$ and may indeed have many different
irreducible components with different dimensions (cf. \cite{KimPandharipande}). 
Nonetheless, it is a
remarkable fact  \cite{Behrend, BehrendFantechi, LiTian}
that $\overline{{\mathcal M}}_{g,n}(X,\beta)$ has a
\lq virtual fundamental class' $[\overline{{\mathcal M}}_{g,n}(X,\beta)]^{\mbox{vir}}$
lying in the expected dimension
$$3g-3 + n + (1-g)\dim X + \int_{\beta} c_1(TX)$$
of $\overline{{\mathcal M}}_{g,n}(X,\beta)$. Gromov-Witten
invariants (originally developed mainly in the case $g=0$ when
$\overline{{\mathcal M}}_{g,n}(X,\beta)$ is more tractable, but now
also studied when $g>0$) are obtained by evaluating cohomology classes
on $\overline{{\mathcal M}}_{g,n}(X,\beta)$ against this virtual
fundamental class. 

The cohomology classes used are of two types. Recall that if $1 \leq j \leq n$
the Witten
class $\psi_j \in H^2(\overline{{\mathcal M}}_{g,n};\QQ)$ is the first
Chern class of $s_j^*(\omega_{g,n})$, where $s_j$ is the $j$th tautological
section of the forgetful map from
$\overline{{\mathcal M}}_{g,n+1}$ to $\overline{{\mathcal M}}_{g,n}$
and $\omega_{g,n}$ is the relative dualising sheaf of this forgetful
map. In a similar way, using the forgetful map from
$\overline{{\mathcal M}}_{g,n+1}(X,\beta)$ to
$\overline{{\mathcal M}}_{g,n}(X,\beta)$, we can define $\Psi_j \in H^2(\overline{{\mathcal M}}_{g,n}(X,\beta);\QQ)$ (and $\Psi_j$ is not quite the pullback of $\psi_j$
via the forgetful map $\pi: \overline{{\mathcal M}}_{g,n}(X,\beta) \to
\overline{{\mathcal M}}_{g,n}$ because of the collapsing process
in the definition of $\pi$). We can also pull back cohomology
classes on $X$ via the evaluation maps $ev_j:\overline{{\mathcal M}}_{g,n}(X,\beta)
\to X$ which send a stable map $f:\Sigma \to X$ to the image $f(p_j)$ of the
$j$th marked point $p_j$ of $\Sigma$ for $1 \leq j \leq n$. 

Gromov-Witten invariants for $X$ are given by integrals
$$\int_{[\overline{{\mathcal M}}_{g,n}(X,\beta)]^{\mbox{vir}}} ev_1^*(\alpha_1) \ldots
ev_n^*(\alpha_n)$$
of classes of the second type $ev_j^*(\alpha_j)$, where $\alpha_1,\ldots,\alpha_n \in
H^*(X;\QQ)$, against the virtual fundamental class of $\overline{{\mathcal M}}_{g,n}(X,\beta)$,
while descendent Gromov-Witten invariants are of the form
$$\int_{[\overline{{\mathcal M}}_{g,n}(X,\beta)]^{\mbox{vir}}} \Psi_1^{k_1}\ldots \Psi_n^{k_n} ev_1^*(\alpha_1) \ldots
ev_n^*(\alpha_n)$$
for nonnegative integers $k_1,\ldots, k_n$, not all zero. More generally, instead of
integrating against $[\overline{{\mathcal M}}_{g,n}(X,\beta)]^{\mbox{vir}}$ to
get rational numbers one can consider the image in $H^*(\overline{{\mathcal M}}_{g,n};\QQ)$
 of the product $\Psi_1^{k_1}\ldots \Psi_n^{k_n} ev_1^*(\alpha_1) \ldots
ev_n^*(\alpha_n)$ under the virtual
pushforward map associated to $\pi: \overline{{\mathcal M}}_{g,n}(X,\beta) \to
\overline{{\mathcal M}}_{g,n}$.

When $X$ is a single point, the descendent Gromov-Witten invariants reduce to the
integrals
$$\int_{\overline{{\mathcal M}}_{g,n}} \psi_1^{k_1} \ldots \psi_n^{k_n}.$$
Witten \cite{Witten} conjectured relations between these integrals (later
proved by Kontsevich \cite{Kontsevich} via a combinatorial description
of $\overline{{\mathcal M}}_{g,n}$) which enable
them to be calculated recursively. Witten's conjecture can be formulated
in terms of the formal power series
$$F_g = \sum_{n\geq 0} \frac{1}{n!} \sum_{k_1,\ldots,k_n \geq 0}
\int_{\overline{{\mathcal M}}_{g,n}}
\psi_1^{k_1} \ldots \psi_n^{k_n} t_{k_1} \ldots t_{k_n}$$
in $\QQ[[t_0,t_1,\ldots]]$: 
it says that $\exp(\sum_{g\geq 0} F_g)$
satisfies a system of differential equations called the Virasoro
relations.

Witten's conjecture has been generalised by Eguchi, Hori and Xiong (with
an extension by Katz)
\cite{CK, EHX, GetzlerVirasoro, GP} to provide relations between Gromov-Witten invariants
and their descendents for general nonsingular projective varieties $X$. Their
generalisation is called the Virasoro conjecture for $X$, since it says that
a certain formal expression (the \lq total Gromov-Witten potential') $Z^X$ in the
Gromov-Witten invariants and their descendents satisfies a system of differential
equations
$${\mathcal L}_k Z^X = 0 \mbox{ for $k \geq -1$}$$
where the differential operators ${\mathcal L}_k$ satisfy the commutation relations
$[{\mathcal L}_k,{\mathcal L}_{\ell}] = (k-\ell){\mathcal L}_{k+\ell}$ and hence span a Lie subalgebra 
of the Virasoro algebra isomorphic to the Lie algebra of polynomial
vector fields in one variable (with ${\mathcal L}_k$ corresponding to
$-x^{k+1}{\rm d}/{\rm d}x$). Dubrovin and Zhang \cite{DZ} have proved
that the Virasoro conjecture determines the Gromov-Witten invariants of $X$
when $X$ is homogeneous.

Getzler and Pandharipande \cite{GP} showed that part of Faber's conjectures
on the structure of the tautological ring of ${\mathcal M}_g$ (the
proportionality formulas) would follow
from the Virasoro conjecture for $X= \CC \PP^2$, and Givental \cite{Givental}
has recently found a proof of the Virasoro conjecture for a class of varieties
which includes all complex projective spaces, thus completing the proof of the
proportionality formulas.

Other methods for finding relations between Gromov-Witten invariants include
the Toda conjecture \cite{GetzlerToda,Getzler2002,OkounkovPandharipande,OP} and
exploitation of intersection theory on $\overline{{\mathcal M}}_{g,n}$
and localisation methods \cite{Bertram2000,Bertram??,
CK,EdidinGraham,FaberPandharipande1999,FP,GraberPandharipande,GV,Kontsevich1995}, which have
been very powerful in enumerative geometry.

\renorm
\section{Moduli spaces of bundles over curves}

\newcommand{\mnd}{{{\mathcal B}_{\Sigma}(r,d)}}

Another very well studied family of moduli spaces is given
by the moduli spaces $\mnd$ of stable holomorphic vector
bundles $E$ of rank $r$ and degree $d$ over a fixed nonsingular
complex projective curve $\Sigma$ of genus $g \geq 2$. 
When $r$ and $d$ are coprime $\mnd$ is a nonsingular complex
projective variety; when $r$ and $d$ have a common factor
then $\mnd$ is nonsingular but not projective, and it has
a natural compactification $\overline{{\mathcal B}}_{\Sigma}(r,d)$
which is projective but singular (except when $g=r=2$) \cite{G2,N2}.
If the curve $\Sigma$ is allowed to
vary as well as the bundle $E$ over $\Sigma$ then we obtain a
\lq universal' moduli space of bundles ${\mathcal B}_g(r,d)$, which maps to the moduli space
${\mathcal M}_g$ of nonsingular curves of genus
$g$ with fibre $\mnd$ over $[\Sigma]$. Pandharipande \cite{Pandharipande}
has shown that ${\mathcal B}_g(r,d)$ has a compactification
$\overline{{\mathcal B}}_g(r,d)$ which maps to $\overline{{\mathcal M}}_g$
with the fibre over $[\Sigma] \in {\mathcal M}_g$ given by
$\overline{{\mathcal B}}_{\Sigma}(r,d)$.

In the case when $r$ and $d$ are coprime we have a good understanding
of the structure of the cohomology ring $H^*(\mnd;\ZZ)$, and this
understanding is particularly thorough 
when $r=2$ \cite{B,KN,ST,Z}. For arbitrary $r$ it is known that the cohomology has no
torsion \cite{AB} and inductive formulas \cite{AB,DR,HN} 
as well as explicit formulas \cite{delBano,LaumonR}  
for computing the Betti
numbers are available.
There is a simple set of
generators for the cohomology ring \cite{AB} and there are explicit formulas for the
intersection pairings between polynomial expressions in these generators,
which in principle determine all the relations by Poincar\'{e} duality
\cite{D,JK2,T}. There is also an elegant description of a complete
set of relations among the generators when $r=2$
\cite{B,KN,ST,Z}, partially motivated by a conjecture of Mumford \cite{K},
and there is a generalisation when $r>2$ which is somewhat less
elegant \cite{EK}. 

When $r$ and $d$ are not coprime the structure
of the cohomology ring $H^*(\mnd;\ZZ)$ is a little more difficult to
describe; for example, the induced Torelli group action 
on $H^*(\mnd;\QQ)$ is nontrivial \cite{CLM}, whereas when $r$ and $d$ are coprime
the Torelli action is trivial and the mapping class group acts via
representations of $Sp(2g;Z)$ which are easy to determine. However even in
this case information
is available on the intersection cohomology of the compactification $\overline{{\mathcal B}}_{\Sigma}
(r,d)$ of ${{\mathcal B}}_{\Sigma}(r,d)$ and the cohomology of another compactification $\tilde{{\mathcal B}}_{\Sigma}(r,d)$ of ${{\mathcal B}}_{\Sigma}(r,d)$ with only orbifold
singularities: for example, there are formulas for the Betti numbers in both cases
\cite{K5} 
and their intersection pairings \cite{JKKW,Kiem}, and the mapping class group again
acts via representations of $Sp(2g;\ZZ)$ \cite{Nelson}.

One of the main reasons for our good understanding of the moduli
spaces $\mnd$ (and their compactifications  $\overline{{\mathcal B}}_{\Sigma}
(r,d)$ and $\tilde{{\mathcal B}}_{\Sigma}(r,d)$ when $r$ and $d$ have a common
factor) is that they can be constructed as quotients, in the sense
of geometric invariant theory \cite{GIT}, of well behaved spaces whose properties
are relatively easy to understand. Similar techniques could in principle
be used to study the moduli spaces of stable curves $\overline{{\mathcal M}}_{g}$
and $\overline{{\mathcal M}}_{g,n}$, as well as Pandharipande's compactification
$\overline{{\mathcal B}}_g(r,d)$ of the universal moduli 
space of bundles ${{\mathcal B}}_g(r,d)$, since they too can be
constructed using geometric invariant theory. In practice this has
not succeeded except in very special cases because, in contrast to
the case of $\mnd$, we do not have quotients of well behaved spaces which are easy to analyse. 
However as our understanding of the moduli spaces $\overline{{\mathcal M}}_{g,n}(X,\beta)$ of stable maps becomes increasingly well developed, 
and in particular localisation techniques are used with greater and greater effect,
perhaps the techniques available for studying the cohomology of geometric invariant theoretic
quotients will provide an additional approach to the cohomology of the moduli
spaces $\overline{{\mathcal M}}_{g}$
and $\overline{{\mathcal M}}_{g,n}$ which can be added to the plethora of methods
already available.

\end{document}